\documentclass[11pt,a4paper]{article}
\usepackage[T1]{fontenc}
\usepackage[utf8]{inputenc}
\usepackage{lmodern}
\usepackage{amsmath,amssymb,amsthm,mathtools}
\usepackage{geometry}
\usepackage{microtype}
\usepackage{enumitem}
\usepackage{fancyhdr}
\usepackage{needspace}
\usepackage{array}
\usepackage{float}
\usepackage[hidelinks]{hyperref}
\usepackage{orcidlink}
\numberwithin{equation}{section}
\newtheorem{theorem}{Theorem}[section]
\newtheorem{lemma}[theorem]{Lemma}
\newtheorem{proposition}[theorem]{Proposition}
\newtheorem{corollary}[theorem]{Corollary}

\newcommand{\rad}{\operatorname{rad}}

\DeclareRobustCommand{\Udots}{%
  \vcenter{\offinterlineskip
    \halign{%
      \hbox to .75em{##}\cr
      \hfil.\cr\noalign{\kern.15ex}
      \hfil.\hfil\cr\noalign{\kern.15ex}
      .\hfil\cr\noalign{\kern.25ex}}%
  }%
}

\hypersetup{pdftitle={Nonconstant Congruence Speeds of Integer Tetration in Non-Squarefree Radices},pdfauthor={Marco Ripà},pdfsubject={Integer tetration, congruence speeds, non-squarefree radices}}
\title{Nonconstant Congruence Speeds of Integer Tetration in Non-Squarefree Radices}
\author{Marco Rip\`a~\orcidlink{0000-0002-6036-5541}\\[4pt]
\small Independent Researcher, Rome, Italy}
\date{November 4, 2026}
\begin{document}
\maketitle
\begin{abstract}
We further develop the theory of congruence speed for integer tetration in non-squarefree numeral systems, building on Appendix~B of [M. Rip\`a and G. Di Pietro (2026), A compact notation for peculiar properties characterizing integer tetration, Notes on Number Theory and Discrete Mathematics, 32(3), 477--496]. For bases coprime to the radix, we establish an eventual floor-difference representation and determine the exact spectrum of attainable least periods. We construct explicit arithmetic progressions containing infinitely many prime bases with nonconstant congruence speed for every non-squarefree radix, thereby proving the coprime-infinitude conjecture posed in Appendix~B. We also derive a uniform quadratic bound for the least coprime base with nonconstant congruence speed and establish exact natural densities of nonstabilizing bases, both among integers coprime to the radix and among those not divisible by its radical, together with the corresponding relative density among primes. Finally, we give necessary and sufficient conditions for the perfect-power congruence-speed identity (Equation (2.12) of the cited paper) to hold in arbitrary radices, determine the exact least eventual period when stabilization fails, and establish a separate criterion for the boundary case of the original parameter restriction.
\end{abstract}

\section{Introduction}
The concept of congruence speed was introduced in 2011 \cite{Book2011} to describe the rate at which the rightmost digits of integer tetration become stable as the height of the corresponding power tower increases. Initially investigated in the decimal numeral system, this property was subsequently studied in greater detail in \cite{Earlier2020,EarlierFormula,StableDigits}, leading to explicit formulas for the constant congruence speed of integer tetration bases and several related results. In particular, \cite{EarlierFormula} proved the existence of infinitely many prime bases having any prescribed positive constant decimal congruence speed, making use of Dirichlet's theorem on primes in arithmetic progressions \cite{Selberg1949,Shapiro1950}.

A natural question arising from these investigations is whether the congruence speed necessarily stabilizes when the decimal radix is replaced by an arbitrary integer $r>1$. This problem was addressed more generally by Rip\`a and Di Pietro \cite{RipaDiPietro}, whose Appendix~B extended the theory to squarefree numeral systems and established that squarefreeness characterizes the universal stabilization of congruence speeds among bases not divisible by the radical of the radix. However, the behavior in non-squarefree systems remained only partially understood, and the existence of infinitely many coprime bases with nonconstant congruence speed for every such radix was explicitly left open as Conjecture~3.1.

In the present paper, we investigate the non-squarefree case by applying lifting-the-exponent (LTE) techniques to the successive differences of integer tetration towers. Rather than restricting our attention to the existence of exceptional bases, we determine the eventual structure of their congruence speeds, showing that these follow periodic patterns whose least periods can be characterized exactly in terms of the prime factorization of the radix. We then establish the aforementioned conjecture by constructing explicit arithmetic progressions containing infinitely many prime bases with nonconstant congruence speed, and obtain an upper bound for the smallest such coprime base.

We also derive exact natural densities of bases with nonconstant congruence speed, both among those coprime to the radix and among the larger class not divisible by its radical. These results are complemented by a corresponding density theorem for prime bases, using the prime number theorem for arithmetic progressions \cite{Davenport}. 
Finally, returning to another result of \cite{RipaDiPietro}, we obtain necessary and sufficient conditions for its perfect-power congruence-speed identity to hold in arbitrary radices, determine the eventual periodic behavior when stabilization fails, and establish a separate criterion for the boundary case of the original parameter restriction.

\section{Definitions and notation}
Throughout this paper, we follow the notation for integer tetration and congruence speed adopted in \cite{RipaDiPietro}, introducing additional definitions as needed.

Let $a,r\ge2$ and $b\ge0$ be integers. We denote by ${}^{b}a$ the integer tetration of base $a$ and height $b$, with ${}^{0}a=1$, ${}^{1}a=a$, and ${}^{b+1}a=a^{{}^{b}a}$ for $b\ge1$. For any prime $p$ and positive integer $n$, let $\nu_p(n)$ be the $p$-adic valuation of $n$, namely the largest nonnegative integer $j$ such that $p^j\mid n$.

The number of stable digits between consecutive tetration heights is given by
\begin{equation}\label{eq:sh}
 s_b^{[r]}(a):=
 \max\{k\ge0:{}^{b+1}a\equiv{}^{b}a\pmod{r^k}\}
 =\min_{p\mid r}
 \left\lfloor
 \frac{\nu_p({}^{b+1}a-{}^{b}a)}{\nu_p(r)}
 \right\rfloor,
\end{equation}
where $p$ ranges over the prime divisors of $r$. The congruence speed at height $b\ge1$ is then defined by
\begin{equation}\label{eq:vh}
 v_b^{[r]}(a) := s_b^{[r]}(a)-s_{b-1}^{[r]}(a).
\end{equation}
We say that the congruence speed stabilizes if $v_b^{[r]}(a)$ is eventually constant. In that case, $\bar b$ denotes the least height from which the congruence speed remains constant, as in \cite{RipaDiPietro}.

Unless otherwise stated, we assume $\gcd(a,r)=1$.

For every prime $p\mid r$, let $\operatorname{ord}_p(a)$ denote the multiplicative order of $a$ modulo $p$, and define
\begin{equation}\label{eq:cp}
c_p(a):=
\begin{cases}
\nu_p\bigl(a^{\operatorname{ord}_p(a)}-1\bigr) & \text{if }p\text{ is odd}\\
\nu_2(a^2-1)-1 & \text{if }p=2
\end{cases}.
\end{equation}
These constants agree with the prime-radix constant congruence speeds given in (3.4) of \cite[Appendix~B]{RipaDiPietro}. For odd primes $p$, LTE also gives $c_p(a)=\nu_p(a^{p-1}-1)$.

The rational LTE slope is defined by
\begin{equation}\label{eq:lambda}
 \lambda_r(a):=\min_{p\mid r}\frac{c_p(a)}{\nu_p(r)}.
\end{equation}
We also write $\rad(r):=\prod_{p\mid r}p$ for the radical of $r$.

Whenever the following set is nonempty, we define
\begin{equation}\label{eq:B}
 B(r):=\min\{a\in\mathbb Z:a\ge2,\ \gcd(a,r)=1,\ \lambda_r(a)\notin\mathbb Z\},
\end{equation}
which gives the least coprime base with nonconstant congruence speed. As $r$ ranges over the non-squarefree integers listed in OEIS A013929, the corresponding values $B(r)$ form OEIS A390576. The related sequence A390535 also allows bases not coprime to $r$ and therefore requires the extended LTE criterion, in which primes dividing $a$ are omitted from the minimum (and an empty minimum corresponds to nonstabilization). The theorems below concern the coprime setting, except where explicitly noted.

Following the publication of \cite{RipaDiPietro} on October 15, 2026, Alper Ferudun presented an LTE criterion on October 20, 2026, subsequently included in OEIS A390535 and A390576 \cite{OEIS535,OEIS576}. For bases coprime to the radix, this criterion establishes that the congruence speed stabilizes if and only if $\lambda_r(a)$ is an integer; the subsequent theorems build on this criterion and the previously published results.

\section{Eventual periodicity and period spectra}
\subsection{The eventual period}
We begin with a classical modular stabilization result for tetration towers \cite[Corollary~2.11]{Shapiro2007}, providing the basis for the subsequent application of LTE.
\begin{lemma}\label{lem:mod}
For every integer $a\ge2$ and every positive integer $n$, the tower sequence ${}^{b}a$ stabilizes modulo $n$. In particular, there exists $b_0$ such that $n\mid {}^{b+1}a-{}^{b}a$ for all $b\ge b_0$.
\end{lemma}
\begin{proof}
We proceed by strong induction on $n$. The case $n=1$ is immediate. Suppose $n>1$ and consider any prime-power divisor $q^k$ of $n$, where $q$ is prime and $k\ge1$.

If $q\mid a$, then $\nu_q({}^{b+1}a)={}^{b}a\,\nu_q(a)\to\infty$ as $b\to\infty$, so ${}^{b}a$ is eventually zero modulo $q^k$.

If $q\nmid a$, Euler's theorem implies that $a^t\pmod{q^k}$ depends only on $t\pmod{\varphi(q^k)}$. Since $\varphi(q^k)<q^k\le n$, the induction hypothesis shows that ${}^{b-1}a$ is eventually constant modulo $\varphi(q^k)$, and therefore ${}^{b}a=a^{{}^{b-1}a}$ is eventually constant modulo $q^k$.

Finally, applying the Chinese remainder theorem to the prime-power factors of $n$ completes the induction.
\end{proof}
\begin{lemma}[Exact $p$-adic increment]\label{lem:lte}
If $p$ is a prime not dividing $a\ge2$, then for all sufficiently large $b$,
\begin{equation}\label{eq:increment}
\nu_p({}^{b+1}a-{}^{b}a)=
\nu_p({}^{b}a-{}^{b-1}a)+c_p(a).
\end{equation}
\end{lemma}
\begin{proof}
Set $D_b:={}^{b+1}a-{}^{b}a$. Then
$D_b=a^{{}^{b-1}a}(a^{D_{b-1}}-1)$ for $b\ge2$, so
$\nu_p(D_b)=\nu_p(a^{D_{b-1}}-1)$.
For odd $p$, Lemma~\ref{lem:mod} ensures that $\operatorname{ord}_p(a)\mid D_{b-1}$ for sufficiently large $b$.
Since $\operatorname{ord}_p(a)\mid p-1$, this order is not divisible by $p$. The LTE lemma therefore gives
\[
\nu_p(a^{D_{b-1}}-1)=\nu_p(a^{\operatorname{ord}_p(a)}-1)+\nu_p(D_{b-1}).
\]
For $p=2$, the base $a$ is odd and $D_{b-1}$ is eventually even, so the $2$-adic LTE formula gives
\[
\nu_2(a^{D_{b-1}}-1)
=\nu_2(a^2-1)-1+\nu_2(D_{b-1}).
\]
These are exactly \eqref{eq:increment}.
\end{proof}
\begin{theorem}[Eventual mechanical word]\label{thm:mechanical}
For every pair of integers $a,r\ge2$ with $\gcd(a,r)=1$, there exist an integer $b_0\ge1$ and $\beta\in\mathbb Q$ such that, for $b\ge b_0$,
\begin{align}
 s_b^{[r]}(a)&=\lfloor b\lambda_r(a)+\beta\rfloor,\label{eq:floor}\\
 v_b^{[r]}(a)&=\lfloor b\lambda_r(a)+\beta\rfloor-
             \lfloor(b-1)\lambda_r(a)+\beta\rfloor.\label{eq:speedfloor}
\end{align}
If $\lambda_r(a)=u/d$ in lowest terms, with $u,d$ positive integers, then the congruence speed has exact least eventual period $d$ and is eventually constant if and only if $d=1$. When $d>1$, its eventual values are $\lfloor u/d\rfloor$ and $\lceil u/d\rceil$, with the larger value occurring exactly $u-d\lfloor u/d\rfloor$ times in each period.
\end{theorem}
\begin{proof}
By Lemma~\ref{lem:lte}, for each $p\mid r$ we have
$\nu_p({}^{b+1}a-{}^{b}a)=c_p(a)b+\gamma_p$ for every sufficiently large $b$, where $\gamma_p\in\mathbb Z$.
Using \eqref{eq:sh} and the identity $\min_i\lfloor x_i\rfloor=\lfloor\min_i x_i\rfloor$, we obtain
\[
 s_b^{[r]}(a)=\left\lfloor
 \min_{p\mid r}\left(\frac{c_p(a)}{\nu_p(r)}b+\frac{\gamma_p}{\nu_p(r)}\right)\right\rfloor.
\]
Terms with slope greater than $\lambda_r(a)$ eventually cease to realize the minimum. Among the terms with that minimal slope, we set $\beta:=\min\{\gamma_p/\nu_p(r):c_p(a)/\nu_p(r)=\lambda_r(a)\}$. Increasing $b_0$ by one if necessary ensures that the same floor formula applies at both $b$ and $b-1$. This proves \eqref{eq:floor} and \eqref{eq:speedfloor}.
Because $\lambda_r(a)=u/d$, the floor differences are periodic with period $d$, take only the values $\lfloor u/d\rfloor$ and $\lceil u/d\rceil$, and contain exactly $u-d\lfloor u/d\rfloor$ occurrences of the larger value in each period. If the congruence speed has an eventual period $D$, summing its values over consecutive periods and using \eqref{eq:floor} shows that $D\lambda_r(a)$ is an integer. Thus $d\mid D$, proving that $d$ is the least eventual period.
\end{proof}
Here, the eventual period refers to the sequence of congruence speed values, not to the cycle of the rightmost nonstable digit. The initial transient is not claimed to obey \eqref{eq:speedfloor}.

\subsection{Classification of attainable periods}
For a fixed $r>1$, let $\mathcal P(r)$ denote the set of least eventual periods obtained as $a\ge2$ ranges over integers coprime to $r$.
\begin{theorem}[Period spectrum]\label{thm:periods}
For every $r>1$,
\begin{equation}\label{eq:periods}
 \boxed{\mathcal P(r)=\{1\}\cup\bigcup_{p\mid r}\{d\ge2:d\mid\nu_p(r)\}.}
\end{equation}
\end{theorem}
\begin{proof}
Theorem~\ref{thm:mechanical} shows that every possible period is the reduced denominator of one of the ratios $c_p(a)/\nu_p(r)$, hence divides an exponent $\nu_p(r)$.

Conversely, fix a prime $p\mid r$ and an integer $d\ge2$ dividing $\nu_p(r)$.
If $p$ is odd, choose
\begin{equation}\label{eq:oddconstruction}
 a=1+\frac{r}{p^{\nu_p(r)-\nu_p(r)/d}}.
\end{equation}
Then $a\equiv1\pmod{q^{\nu_q(r)}}$ for every prime divisor $q\ne p$ of $r$, while $\nu_p(a-1)=\nu_p(r)/d$. In particular, $\gcd(a,r)=1$, $\operatorname{ord}_p(a)=1$, and $c_p(a)=\nu_p(r)/d$. All the other ratios $c_q(a)/\nu_q(r)$ are at least $1$. Thus, $\lambda_r(a)=1/d$ and the period is $d$.

If $p=2$ and $\nu_2(r)\ge3$, choose
\begin{equation}\label{eq:twoconstruction}
 a=1+\frac{r}{2^{\nu_2(r)/d}}.
\end{equation}
Because $\nu_2(r)-\nu_2(r)/d\ge2$, we have $\nu_2(a-1)=\nu_2(r)-\nu_2(r)/d$ and $\nu_2(a+1)=1$; hence, $c_2(a)=\nu_2(r)-\nu_2(r)/d$. Again, all odd-prime contributions are at least $1$, so $\lambda_r(a)=1-1/d$ and the least eventual period is $d$.

The only remaining possibility is $p=2$, $\nu_2(r)=2$, $d=2$. In this case $r/4$ is odd, and we take
\begin{equation}\label{eq:exceptional}
 a=\frac{3r^2}{8}+1.
\end{equation}
Since $r/4$ is odd, $a\equiv7\pmod{16}$, giving $c_2(a)=3$; also $a\equiv1\pmod{q^{2\nu_q(r)}}$ for every odd prime $q\mid r$, so $c_q(a)/\nu_q(r)\ge2$. Consequently $\lambda_r(a)=3/2$, which has period $2$.

Finally, period $1$ is always realized: for $r\ge3$, $a=r-1$ is coprime to $r$ and has $\lambda_r(a)\in\mathbb Z$ (indeed $\lambda_r(a)=1$); for $r=2$, take $a=3$.
\end{proof}
\begin{corollary}
The eventual congruence speed in a squarefree radix is constant for every coprime base. Every non-squarefree radix admits coprime bases with nonconstant periodic congruence speed.
\end{corollary}
\begin{proof}
By \eqref{eq:periods}, $\mathcal P(r)=\{1\}$ if and only if $\nu_p(r)=1$ for every prime $p\mid r$, which is equivalent to $r$ being squarefree.
\end{proof}

For example, $\mathcal P(64)=\{1,2,3,6\}$ and $\mathcal P(144)=\{1,2,4\}$.

\section{Prime bases and explicit bounds}
\subsection{Prime bases with nonconstant congruence speed}
In \cite[Theorem~3]{EarlierFormula}, Dirichlet's theorem was used to construct infinitely many prime bases with any prescribed positive constant decimal congruence speed. We now adapt that approach to the nonconstant case in arbitrary non-squarefree radices.

\begin{theorem}[Infinitely many primes with nonconstant congruence speed]\label{thm:primes}
For every non-squarefree $r>1$, there exist infinitely many prime numbers $\ell\nmid r$ for which the radix-$r$ congruence speed of tetration base $\ell$ is not eventually constant. In fact, there exists an explicit reduced arithmetic progression all of whose terms have nonconstant congruence speed.
\end{theorem}
\begin{proof}
We exhibit a progression $a_0+kM$ with $\gcd(a_0,M)=1$.

If an odd prime $p$ satisfies $\nu_p(r)\ge2$, set
\[a_0:=1+\frac{r}{p^{\nu_p(r)-1}},\qquad M:=r.\]
Then, for every $k\ge0$, the first valuation is $c_p(a_0+kM)=1$ and all other ratios in \eqref{eq:lambda} are at least $1$. Thus, $\lambda_r(a_0+kM)=1/\nu_p(r)$ is nonintegral.

Otherwise, if $\nu_2(r)\ge3$, choose
\[a_0:=1+\frac{r}{2^{\nu_2(r)-2}},\qquad M:=r.\]
Now $c_2(a_0+kM)=2$ and $\lambda_r(a_0+kM)=2/\nu_2(r)$ is nonintegral, while all odd-prime ratios are at least $1$.

Finally, for the remaining nonsquarefree radices, characterized by $\nu_2(r)=2$ and $r/4$ odd and squarefree, we choose
\[
a_0:=\frac{3r^2}{8}+1,\qquad M:=r^2.
\]
For every $k\ge0$, $a_0+kM\equiv7\pmod{16}$ and $a_0+kM\equiv1\pmod{(r/4)^2}$, whence $\lambda_r(a_0+kM)=3/2$. In each of the three cases, $\gcd(a_0,M)=1$, so Dirichlet's theorem guarantees infinitely many primes in the corresponding arithmetic progression.
\end{proof}

Theorem~\ref{thm:primes} establishes the coprime-infinitude assertion of Conjecture~3.1 in \cite{RipaDiPietro}; we do not address the separate multiplicative inequality for composite squarefree radices here.

\subsection{Bounds for the least coprime base with nonconstant congruence speed}
\begin{theorem}[Uniform quadratic bound]\label{thm:bound}
For every non-squarefree integer $r>1$, the minimum $B(r)$ in \eqref{eq:B} exists and satisfies
\begin{equation}\label{eq:quad}
 \boxed{B(r)\le\frac{3r^2}{8}+1.}
\end{equation}
More precisely:
\begin{enumerate}[label=(\roman*)]
\item If an odd prime $p$ satisfies $\nu_p(r)\ge2$, then
\begin{equation}\label{eq:boundodd}
 B(r)\le1+\frac{r}{p^{\nu_p(r)-1}}<r.
\end{equation}
\item If $\nu_2(r)\ge3$, then
\begin{equation}\label{eq:boundtwo}
 B(r)\le1+\frac{r}{2^{\nu_2(r)-2}}<r.
\end{equation}
\item Otherwise, $\nu_2(r)=2$ and $r/4$ is odd and squarefree, and
\begin{equation}\label{eq:boundexception}
 B(r)\le\frac{3r^2}{8}+1.
\end{equation}
\end{enumerate}
In particular,
\begin{equation}\label{eq:necessary}
 \boxed{B(r)>r\ \Rightarrow\ r=2\rad(r).}
\end{equation}
\end{theorem}
\begin{proof}
We use the explicit bases $a_0$ in Theorem~\ref{thm:primes}, with $k=0$; their nonstabilization is established there. In (i) and (ii) the candidate is strictly less than $r$, and in (iii) it equals the stated quadratic bound. Since $r\ge4$, all three casewise bounds are at most $3r^2/8+1$. If $B(r)>r$, neither (i) nor (ii) can hold, so $\nu_2(r)=2$ and $r/4$ is odd and squarefree; hence $\rad(r)=r/2$, giving $r=2\rad(r)$.
\end{proof}

The bound in \eqref{eq:quad} is attained for $r=4$ and $r=12$ (since $B(4)=7$ and $B(12)=55$). The necessary condition \eqref{eq:necessary} was previously suggested empirically in \cite{MOquestion}.

A natural remaining question is whether the uniform quadratic bound in \eqref{eq:quad} can be improved to $B(r)=O(r(\log r)^C)$ for some absolute constant $C>0$.

\section{Densities of nonconstant congruence speeds}
\subsection{Bases coprime to the radix}
For fixed $r>1$, define the relative density among integers coprime to $r$ by
\begin{equation}\label{eq:densitydef}
 \delta(r):=\lim_{X\to\infty}
 \frac{\#\{a\in\mathbb Z:2\le a\le X,\ \gcd(a,r)=1,\ \lambda_r(a)\notin\mathbb Z\}}
 {\#\{a\in\mathbb Z:2\le a\le X,\ \gcd(a,r)=1\}},
\end{equation}
provided that the limit exists.
\begin{theorem}[Density of nonstabilizing coprime bases]\label{thm:density}
The limit in \eqref{eq:densitydef} exists for every $r>1$ and equals
\begin{equation}\label{eq:density}
\boxed{\delta(r)=
\begin{cases}
\dfrac{r-\rad(r)}{r-1} & \textnormal{if }r\textnormal{ is odd}\\[6pt]
\dfrac{(r+1)(r-\rad(r))}{r(r-1)} & \textnormal{if }r\equiv2\pmod4\\[6pt]
\dfrac{r+1-2\rad(r)}{r-1} & \textnormal{if }4\mid r
\end{cases}.}
\end{equation}
In particular, $\delta(r)>0$ precisely when $r$ is not squarefree. If $r/4$ is an odd squarefree integer, then
\begin{equation}\label{eq:rare}
 \boxed{\delta(r)=\frac1{r-1}.}
\end{equation}
\end{theorem}
\begin{proof}
Let $a$ vary uniformly through positive integers coprime to $r$, in the natural-density sense. In this proof, $\mathbb{P}(E)$ denotes the relative natural density of an event $E$ among integers coprime to $r$. For odd $p\mid r$, the distribution of $c_p(a)$ is
\begin{equation}\label{eq:oddprob}
 \mathbb{P}(c_p(a)\ge k)=p^{1-k},\quad k\ge1.
\end{equation}
Indeed, since $\operatorname{ord}_p(a)\mid p-1$ and $p\nmid(p-1)/\operatorname{ord}_p(a)$, the condition $c_p(a)\ge k$ is equivalent to $a^{p-1}\equiv1\pmod{p^k}$. Exactly $p-1$ units modulo $p^k$ satisfy this congruence, out of $\varphi(p^k)=(p-1)p^{k-1}$ units.
For $p=2$ (so $a$ is odd),
\begin{equation}\label{eq:twoprob}
 \mathbb{P}(c_2(a)\ge k)=2^{2-k},\quad k\ge2,
\end{equation}
since $c_2(a)=\nu_2(a^2-1)-1$ and the condition means $a\equiv\pm1\pmod{2^k}$ for $k\ge2$.
Since these conditions are determined by residues modulo powers of distinct primes, the Chinese remainder theorem shows that their relative natural densities multiply.

Then, for each positive integer $j$, the event $\lambda_r(a)=j$ is the difference between the nested events
$\{c_p(a)\ge j\nu_p(r)\text{ for all }p\mid r\}$ and
$\{c_p(a)\ge j\nu_p(r)+1\text{ for all }p\mid r\}$. Thus,
\begin{equation}\label{eq:lambdaint}
 \mathbb{P}(\lambda_r(a)=j)=\prod_{p\mid r}\mathbb{P}(c_p(a)\ge j\nu_p(r))-\prod_{p\mid r}\mathbb{P}(c_p(a)\ge j\nu_p(r)+1).
\end{equation}
For $r$ odd, \eqref{eq:oddprob} gives
\[\mathbb{P}(\lambda_r(a)=j)=\frac{\rad(r)-1}{r^j},\quad j\ge1,\]
so the total stabilizing density is $(\rad(r)-1)/(r-1)$.

For $4\mid r$, \eqref{eq:oddprob}--\eqref{eq:twoprob} give
\[\mathbb{P}(\lambda_r(a)=j)=\frac{2(\rad(r)-1)}{r^j},\quad j\ge1,\]
so the stabilizing density is $2(\rad(r)-1)/(r-1)$.

If $r\equiv2\pmod4$, the $2$-adic factor has exponent $\nu_2(r)=1$ and, since $c_2(a)\ge2$, we treat the case $j=1$ separately. Direct substitution into \eqref{eq:lambdaint} yields
\[
 \mathbb{P}(\lambda_r(a)=1)=\frac{\rad(r)-2}{r},\qquad
 \mathbb{P}(\lambda_r(a)=j)=\frac{2(\rad(r)-1)}{r^j}\quad(j\ge2).
\]
Summing produces stabilizing density
\[\frac{\rad(r)(r+1)-2r}{r(r-1)}.\]
In all three cases, taking the complement yields \eqref{eq:density}.

To justify rigorously the passage from finite congruence calculations to an infinite union over $j$, we note that $\lambda_r(a)\ge j$ implies $c_p(a)\ge j\nu_p(r)$ for every $p\mid r$, a condition modulo finitely many prime powers. The density of this event decays geometrically as $j\to\infty$ by \eqref{eq:oddprob}--\eqref{eq:twoprob}. Thus, the tails are uniformly negligible, and the infinite sum indeed gives the natural density.
The three expressions vanish for squarefree $r$ and are strictly positive for nonsquarefree $r$. If $r/4$ is odd and squarefree, then $\rad(r)=r/2$, giving \eqref{eq:rare}.
\end{proof}
\begin{corollary}[Relative density within the primes]\label{cor:prime-density}
Let $\pi(X)$ be the number of primes not exceeding $X$. For every $r>1$,
\begin{equation}\label{eq:prime-density}
\lim_{X\to\infty}\frac{\#\{\ell\in\mathbb Z:\ell\le X,\ \ell\textnormal{ is prime},\ \ell\nmid r,\ \lambda_r(\ell)\notin\mathbb Z\}}{\pi(X)}=\delta(r).
\end{equation}
\end{corollary}
\begin{proof}
For any fixed finite collection of bounds on the $c_p$, the corresponding events are finite unions of reduced residue classes modulo a suitable power of $r$. For each fixed modulus, the prime number theorem for arithmetic progressions \cite{Davenport} gives the same conditional relative frequencies among primes as among integers coprime to $r$, whereas Dirichlet's theorem alone guarantees infinitude but not the equidistribution required here.

Approximating the event $\lambda_r(a)\in\mathbb Z$ by $\lambda_r(a)\in\{1,\dots,J-1\}$, the omitted event $\lambda_r(a)\ge J$ has upper relative density at most the prime density of an explicit finite union of progressions, which tends geometrically to zero as $J\to\infty$. Passing to the limit, with the finitely many primes dividing $r$ contributing nothing asymptotically, proves \eqref{eq:prime-density}.
\end{proof}
These density formulas concern the frequency of nonstabilizing bases, rather than the magnitude of the least such base $B(r)$.

\subsection{Bases not divisible by the radical}\label{sec:noncoprime-density}
The preceding density theorem concerns only bases coprime to the radix, as in OEIS A390576. We now also admit $a$ sharing some, but not all, of the prime factors of $r$. The resulting problem is directly related to the admissible bases of OEIS A390535.

We consider the relative density
\begin{equation}\label{eq:Deltadef}
\delta_{\mathrm{adm}}(r):=\lim_{X\to\infty}
\frac{\#\{a\in\mathbb Z:2\le a\le X,\ \rad(r)\nmid a,\ v_b^{[r]}(a)\text{ is not eventually constant}\}}
{\#\{a\in\mathbb Z:2\le a\le X,\ \rad(r)\nmid a\}},
\end{equation}
provided this limit exists. The condition $\rad(r)\nmid a$ is crucial: if $\rad(r)\mid a$, all prime-adic valuations of successive tower differences grow faster than linearly, so the behavior belongs to a different, already understood case. For $\rad(r)\nmid a$ there is at least one prime divisor $p$ of $r$ not dividing $a$, and the extended LTE criterion uses only such primes.

Recall that a divisor $d\mid r$ is \emph{unitary} if $\gcd(d,r/d)=1$; this will be written $d\parallel r$.
The function $\delta(d)$ is the coprime nonstabilization density from Theorem~\ref{thm:density}; it vanishes for squarefree $d$.

\begin{theorem}[General admissible-base density]\label{thm:Delta}
For every $r>1$, the limit in \eqref{eq:Deltadef} exists and is rational. Precisely,
\begin{equation}\label{eq:Delta}
\boxed{\displaystyle
\delta_{\mathrm{adm}}(r)=\frac{1}{\rad(r)-1}
\sum_{\substack{d\parallel r\\ d>1}}
\varphi(\rad(d))\,\delta(d).}
\end{equation}
In particular, if $r/4$ is odd and squarefree, then
\begin{equation}\label{eq:Delta4Q}
\boxed{\displaystyle
\delta_{\mathrm{adm}}(r)=\frac{1}{\rad(r)-1}\sum_{d\mid r/4}\frac{\varphi(d)}{4d-1}.}
\end{equation}
\end{theorem}
\begin{proof}
For each admissible $a$, let $S$ be the nonempty set of prime divisors of $r$ that do not divide $a$, and define $d_S:=\prod_{p\in S}p^{\nu_p(r)}$, which is a nontrivial unitary divisor of $r$.
The natural density of integers with exactly this subset $S$ equals
\[
\prod_{p\in S}\Bigl(1-\frac1p\Bigr)
\prod_{p\mid r,\,p\notin S}\frac1p
=\frac{\prod_{p\in S}(p-1)}{\rad(r)}
=\frac{\varphi(\rad(d_S))}{\rad(r)}.
\]
The density of all admissible integers is
$1-1/\rad(r)$. Hence, the conditional weight of the class $S$ is
$\varphi(\rad(d_S))/(\rad(r)-1)$.

For $b\ge2$ and every prime $p\mid r$ with $p\notin S$, factoring the tower difference gives
\[
\nu_p\bigl({}^{b+1}a-{}^{b}a\bigr)={}^{b-1}a\nu_p(a),
\]
since $a^{{}^{b}a-{}^{b-1}a}-1$ is coprime to $p$. These valuations grow faster than linearly, whereas Lemma~\ref{lem:lte} gives eventual affine growth for every $p\in S$. Thus, only primes in $S$ determine the eventual minimum in \eqref{eq:sh}, and the proof of Theorem~\ref{thm:mechanical} applies to that restricted minimum. For the remaining primes, the distribution of the LTE constants $c_p(a)$, conditional on $p\nmid a$, is exactly the same as for a base coprime to $d_S$. The Chinese remainder theorem makes these local conditions independent. As a consequence, within the class $S$, the nonstabilization density is $\delta(d_S)$. Summing the finitely many classes gives \eqref{eq:Delta}.

If $r/4$ is odd and squarefree, a contributing unitary divisor must contain the full factor $4$; those not containing it are squarefree and have $\delta=0$. The contributing divisors are therefore of the form $4d$ with $d\mid r/4$.

Since $d$ is odd, $\rad(4d)=2d$ and $\varphi(\rad(4d))=\varphi(d)$; Theorem~\ref{thm:density} gives $\delta(4d)=1/(4d-1)$. Substituting these identities into \eqref{eq:Delta} proves \eqref{eq:Delta4Q}.
\end{proof}

\begin{corollary}[Nonexistence of a global radix limit]\label{cor:DeltaLimits}
Along the non-squarefree radices,
\[
\boxed{\displaystyle
\liminf_{\substack{r\to\infty\\ r\textnormal{ not squarefree}}}\delta_{\mathrm{adm}}(r)=0,
\qquad
\limsup_{\substack{r\to\infty\\ r\textnormal{ not squarefree}}}\delta_{\mathrm{adm}}(r)=1.}
\]
\end{corollary}
\begin{proof}
For any odd prime $p$, only bases coprime to $p$ are admissible in radix $p^2$. By Theorem~\ref{thm:density}, we have $\delta_{\mathrm{adm}}(p^2)=\delta(p^2)=p/(p+1)$, and therefore
\[
\lim_{\substack{p\to\infty\\ p\text{ prime}}}\delta_{\mathrm{adm}}(p^2)=1.
\]
On the other hand, applying \eqref{eq:Delta4Q} to $r=4p$, we obtain
\[
\delta_{\mathrm{adm}}(4p)=\frac{1}{2p-1}\left(\frac13+\frac{p-1}{4p-1}\right)
=\frac{7p-4}{3(4p-1)(2p-1)}\sim\frac7{24p},
\]
and therefore
\[
\lim_{\substack{p\to\infty\\ p\text{ prime}}}\delta_{\mathrm{adm}}(4p)=0.
\]
Since $0\le\delta_{\mathrm{adm}}(r)\le1$, these two subsequences prove the result.
\end{proof}

Corollary~\ref{cor:DeltaLimits} alone says nothing about the \emph{frequency of radices} with small $\delta_{\mathrm{adm}}(r)$. The following proposition establishes a uniform positive lower bound on a positive-density set of non-squarefree radices, thereby ruling out convergence to zero in natural density.

\begin{proposition}[A positive-density obstruction]\label{prop:nine}
If $9\mid r$, then
\begin{equation}\label{eq:nine}
\boxed{\delta_{\mathrm{adm}}(r)\ge\frac49.}
\end{equation}
Consequently the set of radices $r$ for which $\delta_{\mathrm{adm}}(r)\ge4/9$ has lower natural density at least $1/9$ among positive integers. Restricted to non-squarefree radices, its lower relative density is at least
\begin{equation}\label{eq:nine-relative}
\frac{1/9}{1-6/\pi^2}>0.
\end{equation}
\end{proposition}
\begin{proof}
Every integer $a\equiv2,4,5,7\pmod9$ is coprime to $3$ and satisfies
$\nu_3(a^{\operatorname{ord}_3(a)}-1)=1$. Hence, its contribution to the extended LTE minimum is $1/\nu_3(r)$, which belongs to $(0,1)$ since $9\mid r$. Every other contributing ratio is strictly positive, so the overall minimum is also in $(0,1)$ and cannot be integral. All these four residue classes are admissible since $3\nmid a$ implies $\rad(r)\nmid a$. They constitute a natural density $4/9$ among all integers, and the set of admissible bases has density at most $1$, giving \eqref{eq:nine}. Finally, multiples of $9$ have natural density $1/9$, whereas the non-squarefree integers have natural density $1-6/\pi^2$, proving the asserted lower bound on the relative density.
\end{proof}

Proposition~\ref{prop:nine} also yields a positive lower bound for the mean nonstabilization density over non-squarefree radices. Indeed, multiples of $9$ have natural density $1/9$ and each satisfies $\delta_{\mathrm{adm}}(r)\ge4/9$, whereas the non-squarefree integers have natural density $1-6/\pi^2$. Thus,
\begin{equation}\label{eq:mean-adm}
\liminf_{R\to\infty}\frac{\displaystyle\sum_{\substack{2\le r\le R\\ r\textnormal{ not squarefree}}}\delta_{\mathrm{adm}}(r)}{\#\{r\in\mathbb Z:2\le r\le R,\ r\textnormal{ not squarefree}\}}\ge\frac{4}{81(1-6/\pi^2)}>0.
\end{equation}
Consequently, the mean nonstabilization density over non-squarefree radices does not tend to zero as $R\to\infty$.

For example, $\delta_{\mathrm{adm}}(12)=17/165$, $\delta_{\mathrm{adm}}(18)=61/102$, $\delta_{\mathrm{adm}}(28)=5/117$, and $\delta_{\mathrm{adm}}(36)=137/210$. The first is about $10.30\%$, compared with $\delta(12)=1/11\approx9.09\%$ among coprime bases. The density $\delta_{\mathrm{adm}}(r)$ concerns bases with $\rad(r)\nmid a$; OEIS A390535 also includes bases satisfying $\rad(r)\mid a$ and $r\nmid a$.

\Needspace{12\baselineskip}
\subsection{Numerical examples and computational checks}
Table~\ref{tab:coprime-examples} lists $B(r)$, $\rad(r)$, and $\delta(r)$ for selected radices.

\begin{table}[H]
\centering
\begin{tabular}{|r|r|r|r|}
\hline
$r$ & $B(r)$ & $\rad(r)$ & $\delta(r)$\\
\hline
4 & 7 & 2 & $1/3$\\
8 & 3 & 2 & $5/7$\\
9 & 2 & 3 & $3/4$\\
12 & 55 & 6 & $1/11$\\
18 & 5 & 6 & $38/51$\\
24 & 5 & 6 & $13/23$\\
28 & 215 & 14 & $1/27$\\
36 & 5 & 6 & $5/7$\\
60 & 199 & 30 & $1/59$\\
72 & 5 & 6 & $61/71$\\
132 & 161 & 66 & $1/131$\\
852 & 7687 & 426 & $1/851$\\
\hline
\end{tabular}
\caption{Selected values of $B(r)$, $\rad(r)$, and $\delta(r)$.}
\label{tab:coprime-examples}
\end{table}

Exact-arithmetic computations tested the constructions in Theorems~\ref{thm:periods} and \ref{thm:bound} for all $392$ non-squarefree radices $r\le1000$, including $505$ constructed periods, without discrepancies.

Table~\ref{tab:admissible-density-examples} compares the coprime density $\delta(r)$ with $\delta_{\mathrm{adm}}(r)$, which allows bases not divisible by $\rad(r)$.

\begin{table}[htbp]
\centering
\begin{tabular}{|r|r|r|}
\hline
$r$ & $\delta(r)$ & $\delta_{\mathrm{adm}}(r)$\\
\hline
$4$ & $1/3$ & $1/3$\\
$8$ & $5/7$ & $5/7$\\
$12$ & $1/11$ & $17/165$\\
$24$ & $13/23$ & $297/805$\\
$64$ & $61/63$ & $61/63$\\
$72$ & $61/71$ & $3909/4970$\\
$81$ & $39/40$ & $39/40$\\
$128$ & $125/127$ & $125/127$\\
\hline
\end{tabular}
\caption{Coprime and admissible-base nonstabilization densities.}
\label{tab:admissible-density-examples}
\end{table}

Enumerations of $50000$ integer bases per radix were consistent with these fractions within finite-sample variation. The density theorems themselves follow from the prime-power residue counts and the Chinese remainder theorem, with the $2$-adic distribution treated separately.

\section{Extensions of the perfect-power congruence-speed identity}\label{sec:comparison}
Following \cite{RipaDiPietro}, for any integer $c\ge1$, let $\nu_r(c):=\max\{s\in\mathbb Z:s\ge0,\ r^s\mid c\}$. For composite $r$, this exponent-counting function is not generally a valuation. In the same paper, (2.11) establishes the lower bound $v_b^{[r]}(a)\ge t$ for bases of the form
\[
a=\bigl(z\cdot r^{t+1}+r^{t-\nu_r(c)}\pm1\bigr)^c
\]
under its stated hypotheses, while (2.12) gives the exact constant congruence speed $t$ for squarefree $r\ge3$. The identity was also stated for arbitrary $r>2$, conditional on stabilization, in \cite{PerfectFormula}. We now determine the exact eventual slope for this family in arbitrary radices, together with necessary and sufficient conditions for stabilization and the least eventual period when stabilization fails.

\subsection{Exact slope and stabilization criterion}
\begin{theorem}[Exact slope and congruence speed stabilization criterion for the perfect-power family]\label{thm:perfect-power}
Let $r\ge2$, $c\ge1$, $z\ge0$, and $t\ge\nu_r(c)+2$ be integers. For either choice of sign, consider
\[
 a:=\bigl(z\cdot r^{t+1}+r^{t-\nu_r(c)}\pm1\bigr)^c.
\]
Then $a>1$, $\gcd(a,r)=1$, and
\begin{equation}\label{eq:power-slope}
 \boxed{\displaystyle
 \lambda_r(a)
 =\lambda_r\!\left(\left(r^{t-\nu_r(c)}\pm1\right)^c\right)
 =t+\min_{p\mid r}\frac{\nu_p\!\left(c/r^{\nu_r(c)}\right)}{\nu_p(r)}.}
\end{equation}
The minimum in \eqref{eq:power-slope} belongs to $[0,1)$. In particular, the following statements are equivalent:
\begin{enumerate}[label=(\roman*)]
\item The congruence speed of $a$ is eventually constant.
\item Its eventual constant congruence speed equals $t$.
\item $\rad(r)\nmid c/r^{\nu_r(c)}$, equivalently, some prime $p\mid r$ does not divide $c/r^{\nu_r(c)}$.
\end{enumerate}
When these conditions hold, the exact identity
\begin{equation}\label{eq:power-identity}
 \boxed{\displaystyle
 v_{\bar b}^{[r]}\!\left(\left(z\cdot r^{t+1}+r^{t-\nu_r(c)}\pm1\right)^c\right)
 =v_{\bar b}^{[r]}\!\left(\left(r^{t-\nu_r(c)}\pm1\right)^c\right)=t}
\end{equation}
holds for all $z\ge0$ and both signs. Otherwise, let the minimum in \eqref{eq:power-slope} equal $u/d$ in lowest terms, with $0<u<d$. The eventual congruence speed then has exact least period $d$, takes only the values $t$ and $t+1$, and has exactly $u$ occurrences of $t+1$ in each period. The slope and these period statistics are independent of $z$ and the choice of sign, although the eventual phases need not agree.
\end{theorem}
\begin{proof}
For every prime $p\mid r$,
\begin{equation}\label{eq:power-local-factor}
 \nu_p\!\left(z\cdot r^{t+1}+r^{t-\nu_r(c)}\right)
 =(t-\nu_r(c))\nu_p(r)\ge2,
\end{equation}
since $1+z\cdot r^{\nu_r(c)+1}$ is coprime to $r$. This also gives $a>1$ and $\gcd(a,r)=1$.

For odd $p\mid r$, the integer $p-1$ is even. Whether the sign is positive or negative, LTE applied to the integer inside parentheses, congruent to $\pm1$ modulo $p$, yields
\[
 c_p(a)=\nu_p(a^{p-1}-1)
 =(t-\nu_r(c))\nu_p(r)+\nu_p(c).
\]
Here, $c_p(a)=\nu_p(a^{p-1}-1)$ follows from the definition \eqref{eq:cp} and the standard odd-prime LTE identity.
If $2\mid r$, \eqref{eq:power-local-factor} ensures that the integer inside parentheses is odd and that the $2$-adic valuations of its two neighboring even integers are $1$ and $(t-\nu_r(c))\nu_2(r)$.
The $2$-adic LTE formula therefore gives
\[
 c_2(a)=\nu_2(a^2-1)-1
 =(t-\nu_r(c))\nu_2(r)+\nu_2(c).
\]
Thus, for all $p\mid r$,
\begin{equation}\label{eq:power-local}
 c_p(a)=(t-\nu_r(c))\nu_p(r)+\nu_p(c).
\end{equation}
Consequently,
\[
 \begin{aligned}
 \lambda_r(a)
 &=t-\nu_r(c)+\min_{p\mid r}\frac{\nu_p(c)}{\nu_p(r)}\\
 &=t+\min_{p\mid r}\frac{\nu_p\!\left(c/r^{\nu_r(c)}\right)}{\nu_p(r)}.
 \end{aligned}
\]
This expression is independent of $z$ and the sign, proving both equalities in \eqref{eq:power-slope}. By maximality of $\nu_r(c)$, the integer $c/r^{\nu_r(c)}$ is not divisible by $r$. Hence, at least one prime $p\mid r$ satisfies $\nu_p(c/r^{\nu_r(c)})<\nu_p(r)$, so the minimum in \eqref{eq:power-slope} lies in $[0,1)$.

By Theorem~\ref{thm:mechanical}, the eventual congruence speed is constant precisely when the minimum in \eqref{eq:power-slope} vanishes, in which case it equals $t$. This happens exactly when some prime $p\mid r$ does not divide $c/r^{\nu_r(c)}$, proving the equivalences and \eqref{eq:power-identity}. When the minimum is positive, the period and frequency assertions follow from the same theorem. For instance, with $r=27$, $c=3$, $t=2$, and $z=0$, the two signs yield the distinct eventual period-$3$ patterns $(3,2,2)$ and $(2,2,3)$, respectively, starting at $b=2$.
\end{proof}

\begin{corollary}[Exact scope of the squarefree-radix identity]\label{cor:perfect-power-uniform}
For fixed $r\ge2$, the identity \eqref{eq:power-identity} holds for every $c\ge1$, every admissible $z,t$, and both signs if and only if $r$ is squarefree. For each non-squarefree $r$, there are infinitely many exponents $c$ for which it holds and infinitely many for which it fails, with $t\ge\nu_r(c)+2$ in both cases.
\end{corollary}
\begin{proof}
For squarefree $r$, the condition in Theorem~\ref{thm:perfect-power} is automatic: $\rad(r)=r\nmid c/r^{\nu_r(c)}$. For $r\ge3$, this recovers (2.12) of \cite{RipaDiPietro}.
For arbitrary $r$, taking $c=r^j$ with $j\ge0$ gives $c/r^{\nu_r(c)}=1$, so the identity holds. If $r$ is non-squarefree, taking instead $c=r^j\rad(r)$ gives $c/r^{\nu_r(c)}=\rad(r)$, and \eqref{eq:power-slope} becomes
\begin{equation}\label{eq:power-obstruction}
 \boxed{\displaystyle
 \lambda_r(a)=t+\frac{1}{\max_{p\mid r}\nu_p(r)}\notin\mathbb Z.}
\end{equation}
Thus, the congruence speed is not eventually constant. The two choices provide infinitely many exponents of each kind and exclude all non-squarefree radices from universal validity.
\end{proof}

\subsection{The boundary case}
\begin{proposition}[The boundary $t=\nu_r(c)+1$]\label{prop:power-boundary}
Let $r\ge3$, $c\ge1$, $z\ge0$, and $t=\nu_r(c)+1$ be integers. For either sign, set $a=(z\cdot r^{t+1}+r^{t-\nu_r(c)}\pm1)^c$. For each choice of $z$ and sign, the identity \eqref{eq:power-identity} holds if and only if
\begin{equation}\label{eq:power-boundary}
\boxed{
\begin{cases}
\displaystyle \rad(r)\nmid\frac{c}{r^{\nu_r(c)}} & \textnormal{if }r\not\equiv2\pmod4\\[6pt]
\displaystyle \frac{\rad(r)}2\nmid\frac{c}{r^{\nu_r(c)}} & \textnormal{if }r\equiv2\pmod4
\end{cases}.}
\end{equation}
The criterion is independent of $z$ and the choice of sign. In particular, the inequality $t>\nu_r(c)+1$ in the original published identity can be weakened to $t\ge\nu_r(c)+1$ uniformly for every odd squarefree radix, but not for every even squarefree radix.
\end{proposition}
\begin{proof}
For every odd prime $p\mid r$, LTE gives, as in the preceding theorem,
\[
 \frac{c_p(a)}{\nu_p(r)}
 =\nu_r(c)+1+\frac{\nu_p\!\left(c/r^{\nu_r(c)}\right)}{\nu_p(r)}.
\]
If $4\mid r$, the same equality holds for $p=2$, since the integer $z\cdot r^{t+1}+r^{t-\nu_r(c)}$ is divisible by $4$. Thus, for $r\not\equiv2\pmod4$, we have $\lambda_r(a)=t$ precisely when at least one prime factor of $r$ does not divide $c/r^{\nu_r(c)}$.

If $r\equiv2\pmod4$, the integer $z\cdot r^{t+1}+r$ has $2$-adic valuation exactly $1$. Its sum or difference with $1$ is odd, and the $2$-adic LTE formula gives
\[
 c_2(a)=\nu_2\!\left(z\cdot r^{t+1}+r\pm2\right)+\nu_2(c)
 \ge2+\nu_2(c)>t.
\]
The $2$-adic contribution cannot force the minimum slope to equal $t$; only odd prime factors of $r$ can do so. This gives the second condition in \eqref{eq:power-boundary}. As in Theorem~\ref{thm:perfect-power}, an eventual constant congruence speed equal to $t$ is equivalent to an LTE slope equal to $t$.

If $r$ is odd and squarefree, the first condition is automatic. By contrast, $r=6$, $c=3$, $z=0$, and $t=1$ give $a=(6\pm1)^3$, whose eventual constant congruence speed is $2$, not $1$.

Therefore, the strict bound cannot be universally relaxed for even squarefree radices.
\end{proof}

\subsection{Examples and computational verification}\label{rem:perfect-power-examples}
Table~\ref{tab:perfect-power-examples} gives examples with $t\ge\nu_r(c)+2$, valid for every $z\ge0$ and either sign, with $a$ as in Theorem~\ref{thm:perfect-power}.

\begin{table}[htbp]
\centering
\begin{tabular}{|r|r|r|r|}
\hline
$r$ & $c$ & $\lambda_r(a)-t$ & Least period\\
\hline
$4$ & $1$ & $0$ & $1$\\
$4$ & $2$ & $1/2$ & $2$\\
$12$ & $6$ & $1/2$ & $2$\\
$12$ & $12$ & $0$ & $1$\\
$72$ & $6$ & $1/3$ & $3$\\
$72$ & $72$ & $0$ & $1$\\
$81$ & $9$ & $1/2$ & $2$\\
\hline
\end{tabular}
\caption{Slope offsets ($\lambda_r(a)-t$) and least eventual periods for the perfect-power family.}
\label{tab:perfect-power-examples}
\end{table}

Theorem~\ref{thm:perfect-power} complements the previously established lower bound $v_b^{[r]}(a)\ge t$ for $b\ge2$: when $\lambda_r(a)>t$, it determines the eventual frequencies and least period of the values $t$ and $t+1$. This does not imply that the complete finite-height congruence speed profiles agree for different $z$ or choices of sign. As a computational check, we tested \eqref{eq:power-local} using exact modular arithmetic for all $2\le r\le159$, $1\le c\le80$, $z\in\{0,1,3\}$, $t\in\{\nu_r(c)+2,\nu_r(c)+3\}$, and both signs, finding no discrepancies among $275520$ prime-local valuation equalities.

\section*{AI Declaration}
This manuscript incorporates AI-generated mathematical results developed through extensive interaction with GPT-6 (OpenAI). The author supplied the research framework and previous human results, directed the investigation, and critically reviewed and refined the proposed mathematical arguments and manuscript. The author assumes full responsibility for its contents.

\end{document}